\documentclass{article}

\usepackage{amssymb}
\usepackage[pdftex]{graphicx}

\newtheorem{theorem}{Theorem}

\newtheorem{definition}{Definition}

\title{Model Complete Expansions of the Real Field by Modular Functions and Forms}
\author{Ricardo Bianconi}
\date{}

\begin{document}

\maketitle

\begin{abstract}
We prove a strong form of model completenes for expansions of the field of real numbers by (the real and imaginary parts of) the modular function J, by the modular forms $E_4$ and $E_6$ and quasimodular form $E_2$ defined in the usual fundamental domain, and the restricted sine function and the (unrestricted) exponential function. This is done using ideas of Peterzil and Starchenko's paper \cite{peterzil-starchenko-wp2004} on the uniform definability of $\wp$ function in $\mathbb{R}_{\mathit{an}}$ (and of the modular function $J$). In the conclusion we pose some open problems related to this work.



\end{abstract}

\section{Introduction}

We prove a strong form of model completeness for some expansions of the field of real numbers by modular functions and forms and also by universal families of elliptic functions.

This type of result was first proved for expansions of the field of real numbers by real analytic functions restricted to the product of intervals by Jan Denef and Lou van den Dries in \cite{denef-vddries1988} (this structure is now called $\mathbb{R}_{\mathit{an}}$), by Lou van den Dries for expansions by the restricted sine and exponential functions in \cite{vandendries1988} and by the author for expansions by restricted Weierstrass $\wp$ and abelian functions in \cite{bianconi1991}. These results depend, among other things, on the fact that the functions used in the definition of the structures are analytic at the boundary (or at infinity). A major brakthrough came in 1991 when Alex Wilkie showed the model completeness of the expansion of the field of real numbers by the unrestricted exponential function (which is not analytic at infinity) in \cite{wilkie1996}. Motivated by this result, but using a different technique, Lou van den Dries, Angus Macintyre and David Marker proved in \cite{vandendries-etc-1996} the model completeness of $\mathbb{R}_{\mathit{an,exp}}$, the expansion of $\mathbb{R}_{\mathit{an}}$ by the unrestricted exponential function. These results were used to conclude the o-minimality of such structures.

In 2004, Ya'acov Peterzil and Sergei Starchenko showed in \cite{peterzil-starchenko-wp2004} the remarkable result that the Weierstrass $\wp$ function, as a function of the variables $z$ and $\tau$ (the parameter), is definable in $\mathbb{R}_{\mathit{an,exp}}$ (but not in $\mathbb{R}_{\mathit{an}}$). In their paper they also have shown that the modular function $J$ is also definable in $\mathbb{R}_{\mathit{an}}$. Their results (and some others discussed in the Section \textit{Concluding Remarks}) have led to the question about model completeness of the expansion by such functions. We prove here the (strongly) model completeness of expansions of the field of real numbers by including in the language finitely many function symbols, containing modular functions and forms and the unrestricted exponential function.

The strategy of our proof starts with the ideas of their paper, transforming the problem to a reduct of $\mathbb{R}_{\mathit{an}}$ (via Eisenstein series) and then we use the ideas from \cite{vandendries-etc-1996} to introduce the unrestricted exponential function in order to obtain the model completeness for expansions by modular functions and forms defined in the fundamental domain $\{\tau\in\mathbb{C}:-1/2\leq\mathfrak{Re}(\tau)<1/2$ and $|\tau|\geq 1\}$.

This work is organised as follows. The following section contains the general elimination of quantifies and model completeness results to be used. In the next section we present the main results on the modular function $J$ and on quasimodular forms. In the final section on concluding remarks we comment on the problem of noneffectiveness of our proofs and present some open problems related to the questions treated in this work. We have included in the bibliography the \textit{Mathematical Reviews} number of the references for easy search.

\textbf{Notation.} We denote by $\mathbb{Z}$ the ring of rational integers, $\mathbb{N}$ the set of nonnegative integers, $\mathbb{R}$ the field of real numbers, $\mathbb{C}$ the field of complex numbers; $\mathfrak{Re}(z)$ and $\mathfrak{Im}(z)$ denote the real and imaginary parts of $z$, the letter $i$ denotes $\sqrt{-1}$,  and $\bar{x}$ denotes the tuple $(x_1,\dots,x_n)$ for some unspecified $n$. Other notations are explained in the text.

\section{Strong Model Completeness Results}

The theory $T$ of a structure $M$ is called \textbf{strongly model complete} if for each formula $\varphi(\bar{x})$ there is a quantifier free formula $\psi(\bar{x},\bar{y})$ such that
$$
T\vdash\forall \bar{x}\,[\varphi(\bar{x})\leftrightarrow\exists\bar{y}\,\psi(\bar{x},\bar{y})]
$$
and $T\vdash\forall\bar{x}\forall\bar{y}\forall\bar{z}(\psi(\bar{x},\bar{y})\wedge\psi(\bar{x},\bar{z})\to\bar{y}=\bar{z})$. That is, each formula is equivalent to an existential formula witnessed by a unique element (for each tuple $\bar{x}$).

A set in a structure is called \textbf{strongly definable} if it is defined by such an existential formula.

We first recall a general strong model completeness result for Weierstrass systems with the (unrestricted) exponential functions and then show the particular cases of modular functions and modular forms.

\subsection{Weierstrass Systems}

In this section $K$ stands for the field of real or of complex numbers, $\mathbb{R}$ or $\mathbb{C}$, $K[x_1,\dots,x_n]$ for the ring of polynomials with coefficients in $K$ and the shown variables $x_1$, \dots, $x_n$, $K[[x_1,\dots,x_n]]$ the ring of (formal) power series over $K$ and $K\{x_1,\dots,x_n\}$ its subring of power series converging in some polydisk centred at $\mathbf{0}\in K\sp{n}$. We usually denote $\bar{x}_{(n)}$, or simply $\bar{x}$ when the subscript $n$ is clear from the context, for the $n$-tuple of variables $x_1$, \dots, $x_n$.

\begin{definition}\rm
A \textit{Convergent Weierstrass system, W-system} for short, over $K$ is a family of
rings $W=(K\lfloor x_1,\dots,x_n\rfloor)_{n\in\mathbb{N}}$ satisfying the following conditions, for all $n\in\mathbb{N}$:
\begin{itemize}
\item[W1]
$K[\bar{x}]\subseteq K\lfloor\bar{x}\rfloor \subseteq K[[\bar{x}]]$, $K\lfloor\bar{x}\rfloor$ is closed under the action of the permutation on the variables and if $m>n$, then $K\lfloor\bar{x}_{(n)}\rfloor = K\lfloor\bar{x}_{(m)}\rfloor\cap K[[\bar{x}_{(n)}]]$;
\item[W2]
if $f\in K\lfloor\bar{x}\rfloor$ is invertible in $K[[\bar{x}]]$ then its inverse belongs to $K\lfloor\bar{x}\rfloor$;
\item[W3]
\textit{(Weierstrass Division)}
if $f\in K\lfloor\bar{x}_{(n+1)}\rfloor$ and $f(\mathbf{0},x_{n+1})\in K[[x_{n+1}]]$ is nonzero of order $d\in\mathbb{N}$ (that is, $f(\mathbf{0},x_{n+1})= \sum_{i=d}\sp{\infty}c_ix_{n+1}\sp{i}$, with $c_d\neq 0$) then for all $g\in K\lfloor\bar{x}_{(n+1)}\rfloor$, there are $Q\in K\lfloor\bar{x}_{(n+1)}\rfloor$ and
$R_i\in K\lfloor\bar{x}_{(n)}\rfloor$, $0\leq i\leq d-1$, such that $g(\bar{x}_{(n+1)})=Q(\bar{x}_{(n+1)})f(\bar{x}_{(n+1)})+\sum_{i=0}\sp{d-1}R_i(\bar{x}_{(n)})x_{n+1}\sp{i}$.
\item[W4]
$K\lfloor\bar{x}_{(n)}\rfloor\subseteq K\{\bar{x}_{(n)}\}$ and if the element $f\in K\lfloor\bar{x}_{(n)}\rfloor$ converges on the polydisk $D_{\varepsilon}(0)=\{\bar{x}\in \mathbb{C}\sp{n}:|x_1|<\varepsilon,\dots, |x_n|<\varepsilon\}$ then for each $a\in D_{\varepsilon}(0)\cap K\sp{n}$, the series
$$
f_a(\bar{x}_{(n)})=f(a+\bar{x}_{(n)})=\sum_{\mathbf{i}\in\mathbb{N}}\frac{\partial\sp{\mathbf{i}} f\phantom{n}}{\partial \bar{x}_{(n)}\sp{\mathbf{i}}}(a)\frac{\bar{x}_{(n)}\sp{\mathbf{i}}}{\mathbf{i}!}
$$
also belongs to $K\lfloor{\bar{x}_{(n)}}\rfloor$, where, as usual, if $\mathbf{i}=(i_1,\dots,i_n)\in\mathbb{N}\sp{n}$, then $\mathbf{i}!=i_1!\dots i_n!$, $\bar{x}_{(n)}\sp{\mathbf{i}}=x_1\sp{i_1}\dots x_n\sp{i_n}$ and
$$
\frac{\partial\sp{\mathbf{i}}f\phantom{n}}{\partial \bar{x}_{(n)}\sp{\mathbf{i}}}=\frac{\partial\sp{i_1+\dots+i_n}f}{\partial x_1\sp{i_1}\dots\partial x_n\sp{i_n}}.
$$
\end{itemize}
\end{definition}

Observe that in any W-system, $\mathbb{R}\lfloor\bar{x}_{(0)}\rfloor$ is just the real field $\mathbb{R}$.

Given a W-system $\mathcal{W}=(\mathbb{R}\lfloor \bar{x}_{(n)}\rfloor)_{n\geq 0}$, we define for each $n\geq 1$ the ring $\mathbb{R}_W\{\bar{x}_{(n)}\}$ containing all $f\in\mathbb{R}[[\bar{x}_{(n)}]]$ which converge in a neighbourhood of $I\sp{n}$ and for all $a\in I\sp{n}$ the translated series $f_a$ belongs to $\mathbb{R}\lfloor \bar{x}_{(n)}\rfloor$, where $I$ is the interval $[-1,1]$.

Let $L_{W}$ be the first order language which has the usual symbols for the theory of ordered fields ($+$, $\cdot$, $-$, $<$, ${\,}\sp{-1}$, 0, 1) and also for each $n\geq 1$ and each $f\in\mathbb{R}_W\{\bar{x}_{(n)}\}$ a function symbol $F$.

We make $\mathbb{R}$ into an $L_W$-structure interpreting the symbols of $L_W$ as follows:
\begin{itemize}
\item
the symbols $+$, $\cdot$, $-$, $<$, ${\,}\sp{-1}$, 0, 1 have the usual interpretation, with the covention that $0\sp{-1}=0$;
\item
if $F$ corresponds to the series $f\in\mathbb{R}_W\{\bar{x}_{(n)}\}$, then the interpretation $\tilde{F}$ of $F$ is
$$
\tilde{F}(\bar{x})=
\left\{
\begin{array}{lcl}
f(\bar{x}) & \mbox{if} & \bar{x}\in I\sp{n},\\
0 & & \mbox{otherwise.}
\end{array}
\right.
$$
\end{itemize}

We denote $\mathcal{R}_W$ for such structure.

In \cite[Section 4]{denef-vddries1988}  it was proven a quantifier elimination result for the W-system of all convergent power series. The proof also works for any W-system and the details of such proof can be found in \cite[Section 2]{bianconi1991} (with the appropriate modifications).

\begin{theorem}\label{QE-W-system}
The theory of $\mathcal{R}_W$ admits quantifier elimination in the language $L_W$.
\end{theorem}

The next theorem is based on \cite[Theorem 3.12]{bianconi1991} which gives a criterion for strong model completeness of expansions of the field of the real numbers by restricted analytic functions. The conditions imply that the strongly definable complex analytic functions form a Weierstrass system, permitting the use of the results above. (We say that a complex function is strongly definable in the structure if its real and imaginary parts are strongly definable.)

\begin{theorem}\label{model-completeness-W-system}
Let $\hat{R}=\langle\mathbb{R},\mathit{constants},+,-,\cdot,<,(F_{\lambda})_{\lambda\in\Lambda}\rangle$ be an expansion of the field of real numbers, where for each $\lambda\in\Lambda$, $F_{\lambda}$ is the restriction to a compact polyinterval $D_{\lambda}\subseteq\mathbb{R}\sp{n_{\lambda}}$ of a real analytic function whose domain contains $D_{\lambda}$(and defined as zero outside $D_{\lambda}$), such that there exists a complex analytic function $g_{\lambda}$ defined in a neighbourhood of a polydisk $\Delta_{\lambda}\supseteq D_{\lambda}$ and such that
\begin{enumerate}
\item
$g_{\lambda}$ is strongly definable in $\hat{R}$ and the restriction of $g_{\lambda}$ to $D_{\lambda}$ coincides with $F_{\lambda}$ restricted to the same set;
\item
for each $a\in\Delta_{\lambda}$ there exists a compact polydisk $\Delta$ centred at $a$ and contained in the domain of $g_{\lambda}$, such that all the partial derivatives of the restriction of $g_{\lambda}$ to $\Delta$ are strongly definable in $\hat{R}$.
\end{enumerate}
Under these hypotheses, the theory of $\hat{R}$ is strongly model complete.
\end{theorem}

\subsection{Inclusion of the Exponential Function}

We introduce the full exponential function in order to obtain the main results of this work.

The following theorem is essentially contained in \cite{vandendries-etc-1996}, adapted to the context of W-systems. The authors describe generalised power series fields as models of $\mathbb{R}_{\mathit{an}}$ and use Shoenfield's criterion \cite{shoenfield1971} to prove quantifier elimination. Our proof is the same as in \cite{vandendries-etc-1996}, with the obvious modifications, so we omit it.

\begin{theorem}
Suppose that some restriction of $\exp$ and $\log$ belong to the W-system $\mathcal{W}$. Then
$\hat{R}=\langle \mathbb{R},\,\mathit{constants},+,\cdot,(\cdot)\sp{-1},<,(F_{w})_{w\in W},\exp,\log\rangle$ admits quantifier elimination, where $\exp$ and $\log$ denote the unrestricted exponential and logarithm functions with the convention that $\log x=0$ for $x\leq 0$.
\end{theorem}

Now we can state the main result.

\begin{theorem}
Let $\hat{R}$ be the structure described in Theorem \ref{model-completeness-W-system}. We assume that the functions
$$
e(x)=
\left\{
\begin{array}{lcl}
\exp x & \mbox{if} & 0\leq x\leq 1,\\
0 & & \mbox{otherwise};
\end{array}
\right.
$$
$$
s(x)=
\left\{
\begin{array}{lcl}
\sin x & \mbox{if} & 0\leq x\leq \pi,\\
0 & & \mbox{otherwise},
\end{array}
\right.
$$
have representing function symbols in its language.
The expansion $\hat{R}_{\mathrm{exp}}$ of $\hat{R}$ by the inclusion of the (unrestricted) exponential function ``$\exp$'' is strongly model complete.
\end{theorem}

\textbf{Proof:}
The proof of the analogous to Theorem \ref{model-completeness-W-system} is obtained by interpreting it in another structure like the one in Theorem \ref{QE-W-system}. We can easily adapt the results in \cite[Sections 2, 3 and 4]{vandendries-etc-1996} to our context to prove that the structure $\mathcal{R}_W$ of Theorem \ref{QE-W-system} expanded by the unrestricted $\exp$ and $\log$ still admits elimination of quantifiers. Since the relevant functions of the structure $\mathcal{R}_W$ are strongly definable in $\hat{R}$, we obtain the desired result.\hspace*{\fill}$\square$

\medskip

In the following sections we use this result to show the model completeness of expansions of the real field by modular functions and forms.

As a side remark, a byproduct of the analysis of the reference \cite{vandendries-etc-1996} we obtain that for any W-system $W$, the models of the theory of $\mathbb{R}_W$ which contain $\mathbb{R}$ as a subfield are isomorphic to reducts of the models of $\mathbb{R}_{\mathit{an}}$.

\section{Modular Functions and Forms}

We present below some examples of expansions by modular functions and forms, where we use the general results of the previous section.

\subsection{The Modular Function \textit{J}}

We recall some definitions and results about modular functions (see \cite[Chapters 1 and 2]{schoeneberg1974}).

The modular group is the group $SL(2,\mathbb{Z})$ of $2\times 2$ matrices with integer entries and determinant 1 acting on $\mathbb{H}=\{z\in\mathbb{C}:\mathfrak{Im}(z)>0\}$ by M\"obius transformations
$$
\left(
\begin{array}{cc}
a & b \\
c & d
\end{array}
\right)
\cdot z=\frac{az+b}{cz+d}.
$$
The group is generated by the translation $T(z)=z+1$ and inversion $S(z)=-1/z$ and has a standard fundamental domain $\mathfrak{F}=\{z\in\mathbb{H}$: $|z|\geq 1$ and $-1/2\leq \mathfrak{Re}(z)< 1/2\}$.

A holomorphic (meromorphic) function $f:\mathbb{H}\to\mathbb{C}$ is called a modular function if it satisfies $f(\gamma\cdot z)=f(z)$ for all $\gamma\in SL(2,\mathbb{Z})$. The Siegel Modular Function is the function $J:\mathbb{H}\to\mathbb{C}$ defined as the bijective holomorphic map from $\mathfrak{F}$ onto $\mathbb{C}$ given by the Riemann Mapping Theorem with the conditions that $J(i\infty)=\infty$, $J(i)=1$ and $J(\rho)=0$, where $\rho=(-1+i\sqrt{3})/2$, and extended to the whole $\mathbb{H}$ by the application of all $\gamma\in SL(2,\mathbb{Z})$. All modular functions are rational functions of $J$.

Let $J:\tau\in \mathcal{F}\to J(\tau)\in\mathbb{C}$ be the modular function.
We denote $\tilde{J}:q\in\mathbb{C}\sp{\ast}\to\mathbb{C}$ such that $J(\tau)=\tilde{J}(\exp(2\pi i\tau))$. Since we are considering $q=(u+iv)$ as new variables, and not just a symbol representing $\exp(2\pi i\tau)$, we depart from the usual notation (as in \cite{zagier2008}, for instance) distinguishing $J$ as function of the variable $\tau$ and $\tilde{J}$ as function of the variable $q$.

Notice that $\exp(2\pi i\tau)$ maps the closed strip $\{\tau\in\mathbb{C}:|\mathfrak{Re}(\tau)|\leq 1/2$, $\mathfrak{Im}(\tau)\geq\delta\}$ onto de punctured disk $\{q\in\mathbb{C}:0<|q|\leq\exp(-2\pi\delta)\}$.
This function maps the vertical ray $\{z\in\mathbb{H}:\mathfrak{Re}(z)=-1/2$, $\mathfrak{Im}(z)\geq\delta\}$ onto the segment $\{q\in\mathbb{C}:\mathfrak{Im}(q)=0$, $-\exp(-2\pi\delta)\leq\mathfrak{Re}(q)<0\}$, and the ray $\{z\in\mathbb{H}:\mathfrak{Re}(z)=0$, $\mathfrak{Im}(z)\geq \delta\}$ onto the segment $\{q\in\mathbb{C}:\mathfrak{Im}(q)=0$, $0<\mathfrak{Re}(q)\leq\exp(-2\pi\delta)\}$.
We consider the extension of that mapping sending the point $i\infty$ to the value 0, the limit $\lim_{\tau\to i\infty}\exp(2\pi i\tau)$ (the limit is computed inside the strip).

Kurt Mahler showed in \cite{mahler1969} that the $J(\tau)$, $J'(\tau)$ and $J''(\tau)$ are algebraically independent over $\mathbb{C}(\tau)$ and that $J'''(\tau)$ is a rational function of $J(\tau)$, $J'(\tau)$ and $J''(\tau)$ (see \cite[page 450]{mahler1969}). The derivatives of $j$ are not modular functions, but transform as
$J'(\gamma\cdot \tau)=(c\tau+d)\sp{2}J'(\tau)$ and $J''(\gamma\cdot\tau)=(c\tau+d)\sp{4}J''(\tau)+2(c\tau+d)\sp{3}J'(\tau)$.

In the following theorem $\tilde{J}\sp{(0)}=\tilde{J}$, $\tilde{J}\sp{(1)}=\tilde{J}'$ and $\tilde{J}\sp{(2)}=\tilde{J}''$ (the first and second derivatives of $\tilde{J}$ with respect to the complex variable $q$).

\begin{theorem}\label{model-completeness-J}
The following expansions of the field of the real numbers are strongly model complete, where $q=(u+iv)$ and $|q|=\sqrt{x\sp{2}+y\sp{2}}$:
\begin{enumerate}
\item
$\mathcal{R}_{J_1}=\langle\mathbb{R},\mathit{constants},+,-,\cdot,<,(F\sp{(j)}_n,G\sp{(j)}_n)_{n\in\mathbb{N}, j=0,1,2}\rangle$, where  for each $n\in\mathbb{N}$ and $j=0,1,2$,
$$
F\sp{(j)}_n(q)=\left\{
\begin{array}{lcl}
\mathfrak{Re}(\tilde{J}\sp{(j)}(q)) & \mathrm{if} & |q|\leq 1-\frac{1}{n+1};\\
0 &  & \mathrm{otherwise}.
\end{array}\right.
$$
$$
G\sp{(j)}_n(q)=\left\{
\begin{array}{lcl}
\mathfrak{Im}(\tilde{J}\sp{(j)}(q)) & \mathrm{if} & |q|\leq 1-\frac{1}{n+1};\\
0 &  & \mathrm{otherwise}.
\end{array}\right.
$$
\item
$\mathcal{R}_{J_2}=\langle\mathbb{R},\mathit{constants},+,-,\cdot,<$, $\mathrm{exp}\lceil_{[0,1]}$, $\mathrm{log}\lceil_{[1,2]}$, $\mathrm{sin}\lceil_{[-\pi,\pi]}$, $\mathrm{cos}\lceil_{[-\pi,\pi]}$, $\arctan\lceil_{[-1,1]}$, $(\tilde{J}\sp{(j)}_{\mathit{re}},\tilde{J}\sp{(j)}_{\mathit{im}})_{j=0,1,2}\rangle$, where for each $j=0,1,2$,
$$
\tilde{J}\sp{(j)}_{\mathit{re}}(q)=\left\{
\begin{array}{lcl}
\mathfrak{Re}(\tilde{J}\sp{(j)}(q)) & \mathrm{if} & |q|\leq 1-\delta;\\
0 &  & \mathrm{otherwise}.
\end{array}\right.
$$
$$
\tilde{J}\sp{(j)}_{\mathit{im}}(q)=\left\{
\begin{array}{lcl}
\mathfrak{Im}(\tilde{J}\sp{(j)}(q)) & \mathrm{if} & |q|\leq 1-\delta;\\
0 &  & \mathrm{otherwise},
\end{array}\right.
$$
where $\delta=\exp(-\pi\sqrt{3})$ and $f\lceil_{[a,b]}$ means the function $f$ restricted to the interval $[a,b]$ and defined as zero elsewhere, for $f$ being the sine, cossine, exponential and arctangent functions.
\end{enumerate}
\end{theorem}

\textbf{Proof:} Let us firstly remark that the choice of the number $\delta$ in item 2 was made in order to include the image, under the change of variables $q=\exp(2\pi i\tau)$, of the smallest vertical strip containing the fundamental domain of the action of $SL(2,\mathbb{Z})$ on $\mathbb{H}$. We also remark that although $\tilde{J}$ has a single pole at the origin, the following argument applies (changing $q\tilde{J}(q)$ for $\tilde{J}(q)$, when necessary).

\medskip

\textbf{1}. Let $\Delta_n=\{q\in\mathbb{C}:|q|\leq 1-1/(n+1)\}$.

Because the first and second derivatives of the function $\tilde{J}$ are present and the third derivative is a rational function of $q$, $\tilde{J}$, $\tilde{J}'$ and $\tilde{J}''$, the only requirement to be proved is the second condition of Theorem \ref{model-completeness-W-system}, for $f_n$ being each one among $F_n\sp{(j)}$, $G_n\sp{(j)}$, $j=0,1,2$:

\begin{quote}
\textit{2. 
For each $q\in\Delta_n$ there exists a compact polydisk $\Delta$ centred at $q$ and contained in the domain of $f_n$ such that all the partial derivatives of the restriction of $f_n$ to $\Delta$ are strongly definable in $\hat{R}$.}
\end{quote}

We prove only the case of $F_n\sp{(0)}$  and $G_n\sp{(0)}$, because the other cases are similar and left to the reader.

If $q$ belongs to the interior of $\Delta_n$, let $\epsilon=\mathrm{dist}\,(q,\partial \Delta_n)>0$, the distance from $q$ to the boundary of $\Delta_n$, and let $\Delta(q,\epsilon)$ be the closed disk centred at $q$ and radius $\epsilon$. Using $F\sp{(0)}_n$ and $G_n\sp{(0)}$ restricted to $\Delta(q,\epsilon)$ we can strongly define $\tilde{J}$ restricted to such disk.

If $q\in\partial\Delta_n$, a boundary point of $\Delta_n$, then $q$ belongs to the interior of $\Delta_{n+1}$ and the same argument as above applies, using $F_{n+1}\sp{(0)}$ and $G_{n+1}\sp{(0)}$.

\medskip

\textbf{2.} For this item, we only need to show that the extensions of those functions to larger domains are strongly definable. Here we make use of the modularity property of the function $J(\tau)$ and of its derivatives. The change of variables $q=\exp(2\pi i\tau)$ inverts as $\tau=(\log q)/(2\pi i)$, where we choose the principal branch of $\log$ (that is, $\log q>0$ if $q>0$, both real), because if $\tau=iy$, with $y>0$, then $q=\exp(2\pi i\tau)=\exp(-2\pi y)>0$, a real number. If
$$
\tau'=\frac{a\tau+b}{c\tau+d},
$$
then the corresponding $q'=\exp(2\pi i\tau')$ satisfies
$$
q'=\exp\left[2\pi i\left(\frac{a\log q+2\pi ib}{c\log q+2\pi id}\right)\right].
$$

\begin{figure}[!ht]
\begin{center}
\includegraphics[width=0.9\textwidth]{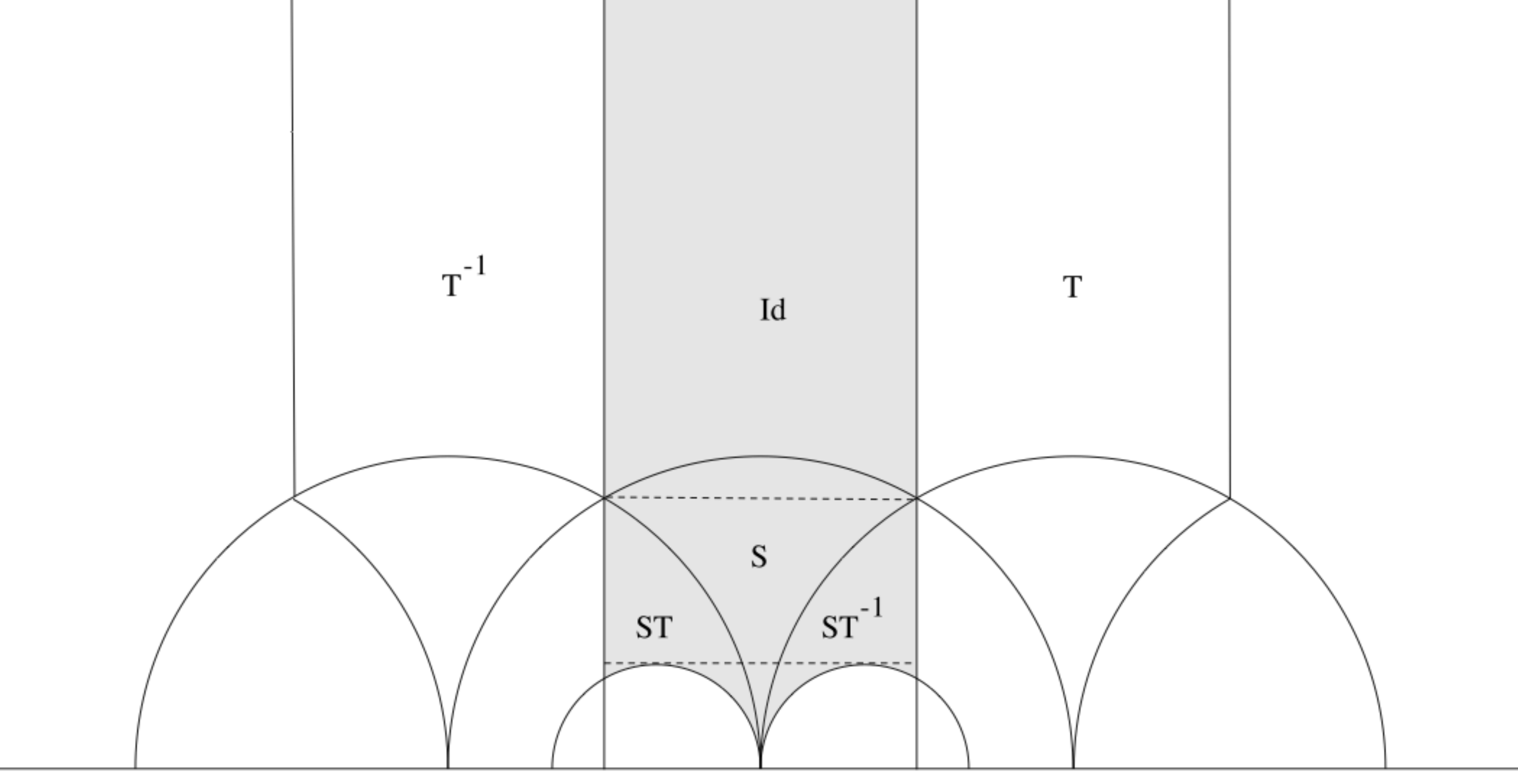}
\end{center}
\caption{Action of $S,ST,ST\sp{-1}\in SL(2,\mathbb{Z})$ on the fundamental domain.}
\label{fig-Action}
\end{figure}

In particular, if we apply this formula using the transformations $\tau_1=S(\tau)=-1/\tau$, $\tau_2=ST(\tau)=-1/(\tau+1)$ and $\tau_3=ST\sp{-1}(\tau)=-1/(\tau-1)$ (the Figure \ref{fig-Action} below helps to explain these choices), we have:
\begin{enumerate}
\item
the image under $S$ of the (closure of the) fundamental domain $\mathfrak{F}$ is the region inside the cilcle $|\tau|=1$, that is $|\tau|\leq 1$ and also outside the cilcles $|\tau\pm 1|=1$;
\item
the image under of $\bar{\mathfrak{F}}$ $ST$ is $\mathfrak{Re}(\tau)\geq -1/2$, $|\tau+1/3|\geq 1/3$ and $|\tau+1|\leq 1$;
\item
the image of $\bar{\mathfrak{F}}$ under $ST\sp{-1}$ is $\mathfrak{Re}(\tau)\leq 1/2$, $|\tau-1/3|\geq 1/3$ and $|\tau-1|\leq 1$.
\end{enumerate}

We conclude that the union of the images of $\bar{\mathfrak{F}}$ under those maps contains the strip $|\mathfrak{Re}(\tau)|\leq 1/2$, $\mathfrak{Im}(\tau)\geq 1/3$.

Notice that the strip $\{\tau\in\mathbb{H}: |\mathfrak{Re}(\tau)|\leq 1/2$, and $\mathfrak{Im}(\tau)\geq i\sqrt{3}/2\}$ is mapped by $\exp(2\pi i\tau)$ onto the punctured disk $\{q\in\mathbb{C}:0<|q|\leq \exp(-\pi\sqrt{3})\}$ and the strip $\{\tau\in\mathbb{H}: |\mathfrak{Re}(\tau)|\leq 1/2$, and $\mathfrak{Im}(\tau)\geq i/3\}$ is mapped by $\exp(2\pi i\tau)$ onto the punctured disk $\{q\in\mathbb{C}:0<|q|\leq \exp(-\pi/3)\}$.

If $\tau\in\bar{\mathfrak{F}}$ and $\tau_1=S(\tau)=-1/\tau$, $\tau_2=ST(\tau)=-1/(\tau+1)$ and $\tau_3=ST\sp{-1}(\tau)=-1/(\tau-1)$, as above, let $q_j=\exp(2\pi i\tau_j)$, $j=1,2,3$, and $q=\exp(2\pi i\tau)$, then we have
$$
q_1=\exp(2\pi i\tau_1)=\exp\left[-2\pi i\left(\frac{2\pi i}{\log q}\right)\right],
$$
$$
q_2=\exp(2\pi i\tau_2)=\exp\left[-2\pi i\left(\frac{2\pi i}{\log q+2\pi i}\right)\right],
$$
$$
q_3=\exp(2\pi i\tau_3)=\exp\left[-2\pi i\left(\frac{2\pi i}{\log q-2\pi i}\right)\right].
$$

\medskip

From the known formulas $\exp(u+iv)=\exp u\,(\cos v+i\sin v)$ and that $\log (u+iv)=\log \sqrt{u\sp{2}+v\sp{2}}+i\arctan(v/u)$, if $u>0$, $\log (u+iv)=\log \sqrt{u\sp{2}+v\sp{2}}+i[\pi/2-\arctan(u/v)]$, if $u<0$ and $v>0$, and $\log (u+iv)=\log \sqrt{u\sp{2}+v\sp{2}}+i[-\pi/2-\arctan(u/v)]$, if $u,v<0$, we see that (the real and imaginary parts of) the points $q_1$, $q_2$ and $q_3$ are strongly definable from the point $q$.
Observe that we only need the restricted exponential and logarithm to define the points $q_j$, $j=1,2,3$, belonging to the set $\{z\in\mathbb{C}:\exp(-\pi/3)\leq|z|<\exp(-\pi\sqrt{3})\}$.

As the derivatives of the functions $\tilde{J}$, $\tilde{J}'$ and $\tilde{J}''$ are strongly definable in that structure, the only thing to be verified is the second item of Theorem \ref{model-completeness-W-system} quoted above, and again the nontrivial step is to consider points at the boundary of the disks.

The functions $\exp$, $\sin$, $\cos$ have addition formulas, and the function $\arctan$ satisfies $\arctan(1/x)=(-1)\sp{\mathrm{sign}(x)}(\pi/2)-\arctan x$, for $x\neq 0$, and $\log ax =\log a+\log x$. These formulas allow us to extend their definitions beyond the boundary points of their restrictions.

The functions $\tilde{J}\sp{(j)}$, $j=0,1,2$, admit the transformation fromulas indicated above (the points $q_1$, $q_2$ and $q_3$, as functions of $q$) in order to extend their definitions beyond the boundary points of their restrictions.

This finishes the proof.\hspace*{\fill}$\square$

\medskip

Now we include of the full exponential and restricted trigonometric functions.

\begin{theorem}\label{model-completeness-J-exp}
Let  $\mathcal{R}_J=\langle\mathbb{R},\mathit{constants},+,-,\cdot,<$, $(J\sp{(k)}_{\mathit{re}},J\sp{(k)}_{\mathit{im}})_{0\leq k\leq 2}$, $\mathrm{sin}\lceil_{[-\pi,\pi]}$, $\mathrm{cos}\lceil_{[-\pi,\pi]}$, $\exp\rangle$, where for each $k=0,1,2$,
$$
J\sp{(k)}_{\mathit{re}}(x,y)=\left\{
\begin{array}{lcl}
\mathfrak{Re}(J\sp{(k)}(x+iy)) & \mathrm{if} & |x|\leq\frac{1}{2},\,y\geq\frac{\sqrt{3}}{2};\\
0 &  & \mathrm{otherwise}.
\end{array}\right.
$$
$$
J\sp{(k)}_{\mathit{im}}(x,y)=\left\{
\begin{array}{lcl}
\mathfrak{Im}(J\sp{(k)}(x+iy)) & \mathrm{if} & |x|\leq\frac{1}{2},\,y\geq\frac{\sqrt{3}}{2};\\
0 &  & \mathrm{otherwise},
\end{array}\right.
$$
and $J\sp{(k)}(x+iy)$ is the $k\sp{th}$ derivative of the (complex) function $J$ with respect to $z$, with $J\sp{(0)}=J$.
The theory of $\mathcal{R}_J$ is strongly model complete.
\end{theorem}

\textbf{Proof:} We interpret the symbols of $\mathcal{R}_{J_2}$ in $\mathcal{R}_{J}$ by strong existential formulas and use the strong model completeness of $\mathcal{R}_{J_2}$.

The change of variables $q=\exp(2\pi i\tau)$, with inverse $\tau=(\log q)/2\pi i$ are strongly definable from each other. The branch of the logarithm function is strongly definable from $\log$ and from $\arctan$, which is itself strongly definable from the restricted sine and cossine.

From these remarks we conclude that the function symbols of the structure $\mathcal{R}_{J_2}$ are strongly definable in $\mathcal{R}_{J}$.

Now, any formula in the language of $\mathcal{R}_J$ can be translated into a formula in the language of $\mathcal{R}_{J_2}$. The strong model completeness of the latter structure implies that such formula is equivalent with a strong existential formula (and this can be done inside $\mathcal{R}_{J}$).  Therefore $\mathcal{R}_J$ is strongly model complete.\hspace*{\fill}$\square$

\subsection{Quasi-modular Forms}

The definitions and claims about modular and quasimodular forms in this section can be found in \cite{zagier2008}.

A \textbf{modular form of (integer) weight} $k>0$ is a holomorphic function $f:\mathbb{H}\to\mathbb{C}$ such that
$$
f\left(\frac{a\tau+b}{c\tau+d}\right)=(c\tau+d)\sp{k}f(\tau)
$$
and of subexponential growth as $\tau\to\infty$ (this implies that $f$ is also holomorphic at $\infty$). We notice that $f(\tau+1)=f(\tau)$, so $f$ admits a \textit{Fourier development}
$\sum_{n\geq 0}a_n\exp(2\pi in\tau)$ (the analycity at $\infty$ implies that the index $n$ varies in the set of nonnegative integers $\mathbb{N}$). It can be shown that these conditions imply that $k\geq 4$.

There are two canonical modular forms of degrees 4 and 6 which generate all the other modular forms:
$$
E_k(\tau)=\frac{1}{2\zeta(k)}\sum_{(m,n)\in\mathbb{Z}\sp{2}\setminus\{(0,0)\}}\frac{1}{(c\tau+d)\sp{k}},
$$
which, for $k>2$, converges absolutely and uniformly in each compact set in $\mathbb{H}$, and $\zeta(k)=\sum_{r\geq 1}(1/r\sp{k})$. The forms 
of even integer weight are polynomials in $E_4$ and $E_6$ with complex coefficients.

For $k=2$ the series does not converge absolutely but we can define
$$
G_2(\tau)=\zeta(2)E_2(\tau)=\frac{1}{2}\sum_{n\neq 0}\frac{1}{n\sp{2}}+\frac{1}{2}\sum_{m\neq 0}\sum_{n\in\mathbb{Z}}\frac{1}{(m\tau+n)\sp{2}},
$$
where the sums cannot be interchanged. This defines a \textbf{quasimodular form of weight 2}, which has the transformation formula
$$
G_2\left(\frac{a\tau+b}{c\tau+d}\right)=(c\tau+d)\sp{2}G_2(\tau)-\pi ic(c\tau+d),
$$
which is not a modular form. Its importance is that the ring $\mathbb{C}[E_2,E_4,E_6]$ contains all the modular forms and is closed under differentiation because of Ramanujan's formulas
$$
E_2'=\frac{E_2\sp{2}-E_4}{12},\ E_4'=\frac{E_2E_4-E_6}{3},\ E_6'=\frac{E_2E_6-E_4\sp{2}}{2},
$$
where $E'_k(\tau)$ denotes $(1/2\pi i)dE_k/d\tau$.

For the representation of these forms in the variable $q$ (Eisenstein series, which can be viewed as Fourier transforms), $\tilde{E}_k(q)$, this derivative corresponds to $qd\tilde{E}/dq$ and so, the ring $\mathbb{C}[q,\tilde{E}_2,\tilde{E}_4,\tilde{E}_6]$ contains the Fourier transforms of all the modular forms and is closed under differentiation:
$$
q\frac{d\tilde{E}_2}{dq}=\frac{\tilde{E}_2\sp{2}-\tilde{E}_4}{12},\ \ q\frac{d\tilde{E}_4}{dq}=\frac{\tilde{E}_2\tilde{E}_4-\tilde{E}_6}{3},\ \ q\frac{d\tilde{E}_6}{dq}=\frac{\tilde{E}_2\tilde{E}_6-\tilde{E}_4\sp{2}}{2}.
$$

As we have done in the case of the modular function $J$ in the last section, we prove here an auxiliary strong model completeness result for the  forms $\tilde{E}_k$ (defined in some disk around the origin) and then, with the help of the exponential function we prove the strong model completeness for the forms $E_k$.

\begin{theorem}
The following expansions of the field of the real numbers are strongly model complete, where $q=(u+iv)$ and $|q|=\sqrt{x\sp{2}+y\sp{2}}$:
\begin{enumerate}
\item
$\mathcal{R}_{M_1}=\langle\mathbb{R},\mathit{constants},+,-,\cdot,<,(F_{k,n})_{n\in\mathbb{N}},(G_{k,n})_{n\in\mathbb{N}, k=2,4,6}\rangle$, where  for each $n\in\mathbb{N}$ and each $k=2,4,6$,
$$
F_{k,n}(q)=\left\{
\begin{array}{lcl}
\mathfrak{Re}(\tilde{E}_k(q)) & \mathrm{if} & |q|\leq 1-\frac{1}{n+1};\\
0 &  & \mathrm{otherwise}.
\end{array}\right.
$$
$$
G_{k,n}(q)=\left\{
\begin{array}{lcl}
\mathfrak{Im}(\tilde{E}_k(q)) & \mathrm{if} & |q|\leq 1-\frac{1}{n+1};\\
0 &  & \mathrm{otherwise}.
\end{array}\right.
$$
\item
$\mathcal{R}_{M_2}=\langle\mathbb{R},\mathit{constants},+,-,\cdot,<$, $\mathrm{exp}\lceil_{[0,1]}$, $\mathrm{log}\lceil_{[1,2]}$, $\mathrm{sin}\lceil_{[-\pi,\pi]}$, $\mathrm{cos}\lceil_{[-\pi,\pi]}$, $\arctan\lceil_{[-1,1]}$, $(\tilde{E}_{\mathit{k,re}},\tilde{E}_{\mathit{k,im}})_{k=2,4,6}\rangle$, where for each $k=2,4,6$,
$$
\tilde{E}_{\mathit{k,re}}(q)=\left\{
\begin{array}{lcl}
\mathfrak{Re}(\tilde{E}_k(q)) & \mathrm{if} & |q|\leq 1-\delta;\\
0 &  & \mathrm{otherwise}.
\end{array}\right.
$$
$$
\tilde{E}_{\mathit{k,im}}(q)=\left\{
\begin{array}{lcl}
\mathfrak{Im}(\tilde{E}_k(q)) & \mathrm{if} & |q|\leq 1-\delta;\\
0 &  & \mathrm{otherwise},
\end{array}\right.
$$
where $\delta=\exp(-\pi\sqrt{3})$ and $f\lceil_{[a,b]}$ means the function $f$ restricted to the interval $[a,b]$ and defined as zero elsewhere, for $f$ being the sine cossine exponential and arctangent functions.
\end{enumerate}
\end{theorem}

\textbf{Proof:} It is similar to the proof of Theorem \ref{model-completeness-J} and left to the reader.\hspace*{\fill}$\square$

\medskip

Now we introduce the full exponential function.

\begin{theorem}
We define $\mathcal{R}_M=\langle\mathbb{R},\mathit{constants},+,-,\cdot,<$, $(E_{\mathit{k,re}}, E_{\mathit{k,im}})_{k=2,4,6}$, $\exp$, $\mathrm{sin}\lceil_{[-\pi,\pi]}$, $\mathrm{cos}\lceil_{[-\pi,\pi]}\rangle$, where for each $k=2,4,6$,
$$
E_{\mathit{k,re}}(x,y)=\left\{
\begin{array}{lcl}
\mathfrak{Re}(E_k(x+iy)) & \mathrm{if} & |\mathfrak{Re}(x)|\leq\frac{1}{2},\,y\geq\frac{\sqrt{3}}{2};\\
0 &  & \mathrm{otherwise}.
\end{array}\right.
$$
$$
E_{\mathit{k,im}}(x,y)=\left\{
\begin{array}{lcl}
\mathfrak{Im}(E_k(x+iy)) & \mathrm{if} & |\mathfrak{Re}(x)|\leq\frac{1}{2},\,y\geq\frac{\sqrt{3}}{2};\\
0 &  & \mathrm{otherwise},
\end{array}\right.
$$
The theory of $\mathcal{R}_M$ is strongly model complete.
\end{theorem}

\textbf{Proof:}
It is similar to the proof of Theorem \ref{model-completeness-J-exp} and is left to the reader.\hspace*{\fill}$\square$

\medskip

Entensions and improvements of the results above are commented in the next section.

\section{Concluding remarks}

All the results in this work rely on the (topological) compactness of the real interval $[-1,1]$ (and of the closed polydisks) which introduces a highly uneffective feature in the proofs.

Recently, Angus Macintyre has shown in \cite{macintyre2003} how to prove an effective version of the model completeness of the expansion of the field of the real numbers by the $\wp$ function with parameter $\tau=i$. His proof can be modified to the cases where $\wp$ admits complex multiplication and perhaps to all cases.
The insight of his work is to deal with the inverse function of $\wp$, which is an elliptic integral and thus it is a Pfaffian function (see details in \cite{macintyre2008}), so it is amenable to the treatment used in his work with Alex Wilkie, \cite{macintyre-wilkie1996}.

On the other hand, his methods do not seem to apply to the modular function $J$ and the modular forms $E_4$ and $E_6$ which are not Pfaffian functions (but are \textit{noetherian} functions) and satisfy nonlinear differential equations involving $E_2$, $E_4$ and $E_6$ in a nontriangularisable way. (P. F. Stiller presents in \cite{stiller1988} second order linear equations satisfied by $E_4\sp{1/4}$ and $E_6\sp{1/6}$ which do no mix the two functions, each reducible to a system of two linear first order equations but, as far as we know, not reducible to the Pfaffian functions type.)

\medskip

These remarks suggest the following \textbf{open problems} related to our work:

\begin{enumerate}
\item
Find proofs of effective model completeness for the expansions of the field of the real numbers by modular functions, forms and elliptic functions (as two variable functions). These proofs would be a big step in the direction of the probable proof of the decidability of their theories.
\item
Find a proof of effective model completeness of expansions by (restricted) noetherian functions (those functions which generate a noetherian differential field over $\mathbb{C}$). The differential algebraic methods of \cite{macintyre-wilkie1996} do not immediately apply here, because they rely on the induction on the length of the \textit{Pfaffian chain} defining the functions and also on effective bounds on the number of connected components of the zero sets of systems of equations with Pfaffian functions. See \cite{khovanskii-gabrielov1998,gabrielov-vorobjov2004} for some work on the complexity of Noetherian Systems.
\item
A useful result needed in the methods of \cite{macintyre-wilkie1996} is an effective bound on the number of connected components of the zero sets of a system of Pfaffian functions. Since the full sine function is a noetherian function, we see the need to restrict their domains to compact (poly)intervals. The measure of the complexity of the noetherian functions have to take into account the size of the polyintervals to which they are restricted.
\item
Prove the same type of results now on the context of $p$-adic numbers. A first step was done by Nathana\"el Mariaule in his 2013 PhD Thesis (supervised by Alex Wilkie), \textit{On the Decidability of the $p$-adic Exponential Ring}, \cite{mariaule2013}, where he proves the effective model completeness of $\mathbb{Z}_p$ expanded by $\exp$ and uses the $p$-adic Schanuel Conjecture to prove the relative decidability of its theory.
\item
The results in \cite{macintyre-wilkie1996}, \cite{macintyre2003,macintyre2008} and \cite{mariaule2013} are relative to the truth some very difficult transcendental theoretic conjectures (Schanuel Conjecture and its modular version). Maybe by a finer analysis of the functions defined by terms on those structures we could find estimates which would circumvent the use of those conjectures.
\item
A question directly related with this work is wether the (full) exponential function is necessary to prove the model completeness of those structures. We believe that, for instance, the full exponential function is not existentially definable from the $\wp$ functions (although it is definable, see \cite[Theorem 5.7, p. 545] {peterzil-starchenko-wp2004}).
\end{enumerate}

\bigskip

\noindent Ricardo Bianconi

\noindent Departamento de Matem\'atica

\noindent Instituto de Matem\'atica e Estat\'istica da Universidade de S\~ao Paulo

\noindent Rua do Mat\~ao, 1010, Cidade Universit\'aria

\noindent CEP 05508-090,
 S\~ao Paulo, SP

\noindent BRAZIL.

\noindent e-mail: bianconi@ime.usp.br

\end{document}